\documentclass[a4paper,12pt,oneside,reqno]{amsart}
\usepackage{amssymb,amsfonts,amsmath,amsthm,amscd, bbm}
\usepackage{graphicx, color}
\usepackage{mathrsfs}
\newcommand{\beq}{\begin{equation}} 
\newcommand{\eeq}{\end{equation}} 
\newcommand{\bea}{\begin{aligned}}
\newcommand{\eea}{\end{aligned}}
\newcommand{\bdm}{\begin{displaymath}}
\newcommand{\edm}{\end{displaymath}}
\newcommand{\barr}{\begin{array}}
\newcommand{\earr}{\end{array}}
\newcommand{\ben}{\begin{enumerate}}
\newcommand{\een}{\end{enumerate}}
\newcommand{\bde}{\begin{description}}
\newcommand{\ede}{\end{description}}

\newcommand{\R}{\mathbb{R}}

\newcommand{\N}{\mathbb{N}}

\newcommand{\PP}{\mathbb{P}}

\newcommand{\EE}{{\sf{E}}}

\newcommand{\defi}{\equiv} 

\newcommand{\be}{\beta}

\newcommand{\la}{\lambda}

\newcommand{\s}{\sigma}

\begin{document}

\title[Solving spin systems]
{Solving spin systems: the Babylonian way.} 

\author[N. Kistler]{Nicola Kistler}
\address{Nicola Kistler \\ J.W. Goethe-Universit\"at Frankfurt, Germany.}
\email{kistler@math.uni-frankfurt.de}

 \date{\today}

\begin{abstract} 
We show that spin systems with generic (ferro- or paramagnetic, or random) interactions are "completely integrable". The approach is worked out, by way of example, for the Sherrington-Kirkpatrick model: we derive an exact formula for the quenched free energy in finite volume which involves an integral over a Gaussian field with correlation structure given by the interaction matrix (with a twist). 
\end{abstract}

\thanks{This work has been supported by a DFG research grant, contract number 2337/1-1. It is my pleasure to thank Stephan Gufler, Goetz Kersting, Markus Petermann, Adrien Schertzer, Marius A. Schmidt,  and Giulia Sebastiani for the endless discussions on the topic of these notes.}

\maketitle

The Sherrington-Kirkpatrick (SK) model \cite{sk} for mean field spin glasses is constructed as follows: for $N \in \N$,  consider centered Gaussians $( g_{ij})_{1 \leq i < j \leq N}$ issued on some probability space $(\Omega, \mathcal F, \PP)$. These Gaussians, {\it the disorder}, are assumed to be all independent and with variance $1/N$. The Ising configuration space $\Sigma_N \defi \{\pm 1\}^N$ is endowed with coin tossing measure, $P_o(\s) \defi 2^{-N}$ for $\s \in \Sigma_N$. We denote by $E_o$  expectation under $P_o$. 

The  SK-Hamiltonian is  
\beq
H_N(\s) \defi \sum_{1 \leq i < j \leq N} g_{ij} \s_i \s_j = \frac{1}{2} \sum_{i,j=1}^N g_{ij} \s_i \s_j\,,
\eeq
with symmetrised random interaction matrix ${\bf G }\defi (g_{ij})$, where $g_{ij} = g_{ji}$, and $g_{ii}=0$.  

The quenched SK-free energy to  inverse temperature $\be>0$ and external field $h \in \R$ is
\beq
N f_N(\be, h) \defi \log E_o \exp\left( \be H_N(\s) + h \sum_{i=1}^N \s_i \right)\,.
\eeq
{\bf Theorem.} {\it  Let ${\sf X} = ({\sf X}_i)_{i=1}^N$ be a centered Gaussian random field  with covariance
\beq \label{covi}
\EE {\sf X}_i {\sf X}_j = 
\begin{cases}
\sum_{k=1}^N \left|g_{ik}\right| \,, & i = j, \\
g_{ij}\,, & i \neq j. 
\end{cases}
\eeq
Then the exact formula for the {quenched}, finite volume SK-free energy holds:
\beq \label{exact_fe}
N f_N(\be, h) = \log \EE \exp\left( \sum_{i=1}^N \log \cosh\left(h + \sqrt{\be} {\sf X}_i \right) \right) - \frac{\be}{2} \sum_{i,j=1}^N |g_{ij}| \,.
\eeq
}
A minor yet key modification of the diagonal therefore turns the Gaussian disorder matrix $\bf G$ into the covariance of yet another Gaussian field.  

Before giving the elementary proof of the theorem, some further comments are in order. 

First we stress that the treatment also yields exact formulas for spin magnetizations, (higher order) correlations etc. under the quenched Gibbs measure. 

We furthermore note that, albeit quenched, formula \eqref{exact_fe} is reminiscent of expressions obtained through replica computations before breaking of the symmetries \cite{mpv}. What is perhaps even more surprising, the approach we implement here relies on a positivity principle combined with a Gaussian integral to decouple the two-body interaction in order to perform the trace over the Ising spins: this is particularly intriguing since the mathematically elusive "replica trick" \cite{kac, sk} relies on  the very same steps,  {\it but executed in different order}, to achieve the very same goal.

We emphasize that neither the Ising nature of the spins, nor the Gaussian character of the disorder play any structural role in the derivation. The method allows, {\it mutatis mutandis}, to derive exact formulas for the Hopfield\footnote{The ensuing formula differs from the one obtained via plain Hubbard-Stratonovich  transformation.} model \cite{hopfield}, the independent set problem \cite{m}, the Edwards-Anderson (EA) model \cite{ea} or, incidentally, also classical systems such as the Ising model on $\mathbb Z_d$ for $d\geq 1$, or, for that matter,  on any (random) graph etc. The approach relies on the {\it Babylonian trick}\footnote{as Res Jost referred to the method of completing the squares.} to decouple the two body interaction, and thus applies in vast generality. (As a matter of fact, the procedure is flexible enough to cover $p$-body interactions \cite{gm} for any $p \geq 2$, or the perceptron \cite{rosenblatt}: how to do this is sketched in the remark at the end of the proof).

The case of the EA-model is particularly interesting, and relevant. The model is to these days deemed intractable, also from the standpoint of theoretical physics: indeed, even numerical simulations have been frighteningly inconclusive so far. Now, a minor variation of formula \eqref{exact_fe} holds true for this model as well: at the risk of being opaque, one simply "perturbs" (some of) the external fields with the boundary conditions (keeping in mind that, in case of the EA-model, the entries of the interaction matrix vanish unless two sites are nearest neighbours). It thus follows that all physically relevant quantities can be efficiently simulated via polynomially many Gaussian random variables only. There is reason to believe that this dramatic reduction of the complexity will help to settle the long-standing debate between the droplet picture of Fisher and Huse \cite{fh} in low dimensions and the replica symmetry breaking scenario of the Parisi theory \cite{mpv}. 

The decoupling is straightforward also in models where the underlying measure ($P_o$) itself is interacting, such as the Heisenberg model \cite{he} for magnetism ($P_o \to$ spherical measure) or, say, the Domb-Joyce model \cite{d} for weakly self avoiding walks ($P_o \to$ law of simple random walk), but the ensuing formulas are more involved.


Given the amount of similarity among the "Babylonian formulas" for the models listed above, it is somewhat puzzling how different choices of the covariance matrix in the underlying Gaussian fields lead to drastically different types of phase transitions. This, if anything, suggests that a rigorous asymptotic evaluation of the Gaussian integrals may not be a simple matter, in general.  One may also expect that different types of covariance structure require radically different tools: when looked through the lenses of the Gaussian representations, classical models (such as Ising) appear as a branch of algebraic spectral theory (the covariance is given by adjacency matrices modified on the diagonal); mean field models (such as SK) are strongly tied to correlated Wigner matrices, whereas realistic spin glasses (such as EA) somehow link with {\it random} algebraic geometry. 
These concerns seem however premature: after all, even basic features of the (random) matrices as in \eqref{covi} have not yet been fully addressed, let alone properties of the associated Gaussian fields. 

Finally, the Babylonian trick stems from an effort to deal directly, and decidedly, with the issue of frustration, ever so substantial yet unwieldy in disordered systems.

\begin{proof}[Proof of the theorem] Some conventions: since the  parameter $N$ is fixed throughout, we will mostly omit it from our notations. In particular we write $\sum_i a_i \defi \sum_{i=1}^N a_i$ and $\sum_{ij} a_{ij} \defi \sum_{i,j=1}^N a_{ij}$ for double sums. We shorten
\beq \label{split_rv}
g_{ij}^+ \defi \max\left\{ g_{ij}, 0\right\}, \qquad g_{ij}^{-} \defi \max\left\{ - g_{ij}, 0\right\}\,,
\eeq
in which case
\begin{itemize}
\item[i)] both $g^+$ and $g^-$ are positive: $g_{ij}^{+} \geq 0$ and $g_{ij}^{-} \geq 0$;
\item[ii)] the decomposition holds true: $ g_{ij} = g_{ij}^{+}- g_{ij}^{-}$;
\item[iii)] similarly, it holds: $ \left| g_{ij }\right| = g_{ij}^{+} + g_{ij}^{-}$.
\end{itemize}
We first split the Hamiltonian
\beq \label{split} \bea
\sum_{i, j} g_{ij} \s_i \s_j & \stackrel{\text{ii)}}{=} \sum_{i, j} g_{ij}^+ \s_i \s_j + \sum_{i, j} \left(-g_{ij}^-\right) \s_i \s_j \\
& = \sum_{i, j} g_{ij}^+ \s_i \s_j + \sum_{i, j} g_{ij}^{-} \left( - \s_i \s_j \right)\,.
\eea \eeq
The {\it Babylonian trick} amounts to writing
\beq \label{baby_one}
\s_i \s_j = \frac{1}{2}(\s_i+ \s_j)^2 - 1, \qquad - \s_i \s_j = \frac{1}{2}(\s_i- \s_j)^2 - 1,
\eeq 
(also using the simplification $\s_i^2= 1$, valid for Ising spins). This implies
\beq \label{split_2}
\eqref{split} = \frac{1}{2} \sum_{i, j} g_{ij}^+ (\s_i+ \s_j)^2 +  \frac{1}{2} \sum_{i, j} g_{ij}^- (\s_i- \s_j)^2 + \boldsymbol g\,,
\eeq
where the $\s$-independent term reads
\beq
{\boldsymbol g} \defi - \sum_{i,j} g_{ij}^+ - \sum_{i,j} g_{ij}^-  = - \sum_{i,j} \left(g_{ij}^+ + g_{ij}^- \right) \stackrel{\text{iii)}}{=} - \sum_{i,j} |g_{ij}|\,.
\eeq
The quenched free energy is thus 
\beq \bea \label{ten}
& N f_N(\be, h) = \log E_o \exp\left(  \frac{\be}{4} \sum_{i,j} g_{ij}^+ (\s_i+ \s_j)^2 +  \frac{\be}{4} \sum_{i, j} g_{ij}^- (\s_i- \s_j)^2 + h \sum_{i} \s_i \right)  + \frac{\be}{2} {\boldsymbol g} \,.
\eea \eeq
We henceforth focus on the quenched partition function 
\beq
Z_N(\be, h) \defi E_o  \exp\left(  \frac{\be}{4} \sum_{i,j} g_{ij}^+ (\s_i+ \s_j)^2 +  \frac{\be}{4} \sum_{i, j} g_{ij}^- (\s_i- \s_j)^2 + h \sum_i \s_i \right) \,.
\eeq
The double sums in the exponential consist of {\it positive} terms: we can thus apply the Hubbard-Stratonovich transformation \cite{h}. To do so we introduce $(\tilde X_{ij}), (\tilde Y_{ij})$ standard Gaussians, all independent, issued on some probability space $(\Omega', \mathcal F', {\sf P})$. We denote by $\EE$ their joint expectation, and shorten
\beq \label{sho}
{\sf X}_{ij} \defi \tilde X_{ij} \sqrt{g_{ij}^+} , \qquad {\sf Y}_{ij} \defi \tilde Y_{ij} \sqrt{g_{ij}^-} \,.
\eeq
We stress that the disorder matrix ${\bf G}$ is quenched. In particular, $g_{ij}$, $g_{ij}^+$ and $g_{ij}^-$ are all constant. Furthermore, the matrices ${\bf G}^+ \defi \left( g_{ij}^+\right),  {\bf G}^- \defi \left( g_{ij}^-\right)$  are by construction symmetric, but the $ {\sf X} \defi ({\sf X}_{ij}),  {\sf Y} \defi ({\sf Y}_{ij})$ are {\it not}. Performing a Hubbard-Stratonovich transformation, and interchanging the order of integration yields
\beq \label{z_int}
Z_N(\be, h) \defi \EE E_o \exp\left(  \sqrt{\frac{\be}{2}} \sum_{i, j} {\sf X}_{ij}  (\s_i+ \s_j) +  \sqrt{\frac{\be}{2}} \sum_{i, j} {\sf Y}_{ij}  (\s_i- \s_j) + h \sum_i \s_i \right) \,.
\eeq
We write the double sums in \eqref{z_int} as 
\beq \bea \label{double_sum}
& \sum_{i, j} {\sf X}_{ij}  (\s_i+ \s_j) +  \sum_{i, j} {\sf Y}_{ij}  (\s_i- \s_j) = \sum_i \left( {\sf X}_{i \bullet} + {\sf X}_{\bullet i} + {\sf Y}_{i \bullet} - {\sf Y}_{\bullet i}  \right) \s_i\,,  
\eea \eeq
with the shorthand notation $A_{ i \bullet} \defi \sum_j a_{ij}$ and $A_{\bullet j} \defi \sum_i a_{ij}$, which we use for any $N\times N$ matrix $A = (a_{ij})$. Using \eqref{double_sum} in \eqref{z_int}, and integrating out the Ising spins yields 
\beq \bea \label{cruciale}
Z_N(\be, h) = \EE \exp\left( \sum_{i=1}^N \log \cosh\left( h+ \sqrt{\frac{\be}{2}} \left( {\sf X}_{i \bullet} + {\sf X}_{\bullet i} + {\sf Y}_{i \bullet} - {\sf Y}_{\bullet i}\right)  \right) \right)\,.
\eea \eeq
This is, in essence, the claim of the theorem, the rest being only a (co)variance check. \\

\noindent The involved random variables are by definition centered:
\beq \label{center}
\EE\left[ {\sf X}_{i \bullet} + {\sf X}_{\bullet i} \right]= 0\,, \, \qquad \EE\left[ {\sf Y}_{i \bullet} - {\sf Y}_{\bullet i} \right] = 0\,.
\eeq
Furthermore, ${\sf X}$ and $ {\sf Y}$ are independent, hence
\beq \label{varxy}
{\sf var}\left[ \left({\sf X}_{i \bullet} + {\sf X}_{\bullet i}\right)  + \left({\sf Y}_{i \bullet} - {\sf Y}_{\bullet i} \right)\right] = {\sf var}\left[ {\sf X}_{i \bullet} + {\sf X}_{\bullet i}\right]  + {\sf var}\left[{\sf Y}_{i \bullet} - {\sf Y}_{\bullet i} \right]\,,
\eeq
By \eqref{center}, the variance equals the second moment: writing out the $\sf X$-contribution we get
\beq \bea
{ \sf var}\left[ {\sf X}_{i \bullet} + {\sf X}_{\bullet i}\right] = \sum_{j,k} \left( \EE {\sf X}_{i j} {\sf X}_{i k}+\EE {\sf X}_{i j} {\sf X}_{ki} + \EE {\sf X}_{ji} {\sf X}_{ik} + \EE {\sf X}_{ji} {\sf X}_{ki}\right)\,.
\eea \eeq
The middle terms vanish (recall also that $g_{ii} = 0$, hence ${\sf X}_{ii} = 0$), whereas first and last are $\neq 0$ only if $k=j$. Using furthermore that $\EE {\sf X}_{ij}^2 = g_{ij}^+ = g_{ji}^+ = \EE {\sf X}_{ji}^2 $ (symmetry of ${\bf G}^+$), we get
\beq \label{2x}
{\sf var}\left[ {\sf X}_{i \bullet} + {\sf X}_{\bullet i}\right] = 2 \sum_j  \EE\left[{\sf X}_{i j}^2\right] \stackrel{\eqref{sho}}{=} 2\sum_j  \EE\left[\left(\sqrt{g_{ij}^+} \tilde X_{ij} \right)^2 \right]  = 2  \sum_j g_{ij}^+ = 2 g_{i\bullet}^+\,.
\eeq
The computation of the second term on the r.h.s. of \eqref{varxy} is just as straightforward: also using, this time, that $\EE {\sf Y}_{ij}^2 = g_{ij}^- = g_{ji}^- = \EE {\sf Y}_{ji}^2$ (symmetry of ${\bf G}^-$), one steadily checks that
\beq \label{2y}
{ \sf var}\left[{\sf Y}_{i \bullet} - {\sf Y}_{\bullet i} \right] = \EE\left[ \left({\sf Y}_{i \bullet} -{\sf Y}_{\bullet i} \right)^2 \right] = 2 \sum_{j} g_{ij}^- = 2 g_{i \bullet}^- \,.
\eeq
Plugging \eqref{2x} and \eqref{2y} in \eqref{varxy} therefore yields
\beq \label{penultimate}
{\sf var}\left[ \left( {\sf X}_{i \bullet} + {\sf X}_{\bullet i}\right)  + \left({\sf Y}_{i \bullet} - {\sf Y}_{\bullet i} \right)\right] = 2 \left( g_{i\bullet}^+ + g_{i\bullet}^- \right) \stackrel{\text{iii)}}{=} 2 \sum_k \left| g_{ik} \right|\,.
\eeq
We next compute, for $i \neq j$, the covariance: 
\beq \bea \label{cov1}
{\sf C}_{ij} & \defi  \EE\left\{  \left( {\sf X}_{i \bullet} + {\sf X}_{\bullet i}\right) + \left({\sf Y}_{i \bullet} - {\sf Y}_{\bullet i} \right)\right\} \left\{ \left({\sf X}_{j \bullet} + {\sf X}_{\bullet j}\right)  + \left({\sf Y}_{j \bullet} - {\sf Y}_{\bullet j} \right) \right\} \\
& =  {\EE\left[ \left({\sf X}_{i \bullet} + {\sf X}_{\bullet i}\right) \left({\sf X}_{j \bullet} + {\sf X}_{\bullet j}\right) \right]} + {\EE\left[ \left({\sf Y}_{i \bullet} - {\sf Y}_{\bullet i} \right) \left({\sf Y}_{j \bullet} - {\sf Y}_{\bullet j} \right) \right] }\,,
\eea \eeq
again since $\sf X$ and $\sf Y$ are independent (and centered). The $\sf X$-contribution reads
\beq \bea \label{first_x}
& \EE\left[ \left( {\sf X}_{i \bullet} +  {\sf X}_{\bullet i}\right) \left( {\sf X}_{j \bullet} +  {\sf X}_{\bullet j}\right) \right] = \EE\left[  {\sf X}_{i \bullet}   {\sf X}_{j \bullet} \right] + \EE\left[ {\sf X}_{i \bullet}  {\sf X}_{\bullet j} \right]+ \EE\left[  {\sf X}_{\bullet i} {\sf X}_{j \bullet}  \right]+ \EE\left[  {\sf X}_{\bullet i}  {\sf X}_{\bullet j} \right].
\eea \eeq
First and last term on the r.h.s. of \eqref{first_x} vanish since for $i \neq j$ the involved random variables are independent and centered, whereas second and third terms read, respectively,
\beq \bea \label{x2}
& \EE\left[ {\sf X}_{i \bullet}  {\sf X}_{\bullet j} \right] = \sum_{k,l} \EE\left[ {\sf X}_{i k}  {\sf X}_{l j} \right] = 
\sum_{k,l} 1_{k=j, l=i} \EE\left[ {\sf X}_{i k}  {\sf X}_{l j} \right]   = \EE\left[ {\sf X}_{i j}^2 \right] = g_{ij}^+, \\
& \EE\left[  {\sf X}_{\bullet i} {\sf X}_{j \bullet}  \right] = \sum_{k,l} \EE\left[  {\sf X}_{k i} {\sf X}_{j l}  \right] = \sum_{k,l}1_{k=j, l=i} \EE\left[  {\sf X}_{k i} {\sf X}_{j l}  \right]  = g_{ji}^+ = g_{ij}^+,
\eea \eeq
again by independence, and symmetry in the last step. Using all this in \eqref{first_x} yields
\beq \label{firstx}
\EE\left[ \left( {\sf X}_{i \bullet} +  {\sf X}_{\bullet i}\right) \left( {\sf X}_{j \bullet} +  {\sf X}_{\bullet j}\right) \right] = 2 g_{ij}^+ \,.
\eeq
The computation of the ${\sf Y}$-contribution in \eqref{cov1} is fully analogous: the upshot reads
\beq \label{ultimo}
\EE\left[  \left({\sf Y}_{i \bullet} - {\sf Y}_{\bullet i} \right) \left({\sf Y}_{j \bullet} - {\sf Y}_{\bullet j} \right) \right] = -2 g_{ij}^- \,.
\eeq
Combining \eqref{firstx} and \eqref{ultimo} yields, for $i \neq j$, 
\beq \bea \label{ultimate}
& {\sf C}_{ij} =  2 g_{ij}^+ - 2 g_{ij}^-  \stackrel{\text{ii)}}{=} 2 g_{ij}.
\eea \eeq
Using \eqref{penultimate} and \eqref{ultimate} in \eqref{cruciale}, together with \eqref{ten}, settles the claim of the theorem. 

\end{proof}

\noindent {\bf Remark.}  { \it The treatment works also for generic $p$-spin interactions. Consider e.g. the  mean field 3-spins (Ising) Hamiltonian
\beq
H_{N, 3}(\s) \defi \sum_{1 \leq i < j < k \leq N} g_{ijk} \s_i \s_j \s_k,
\eeq
with, say, random couplings: one uses the splitting \emph{ii)} and the Babylonian trick
\beq
\s_i \s_j \s_k = \frac{1}{2} \left( \s_i + \s_j \s_k \right)^2- 1, \qquad  - \s_i \s_j \s_k = \frac{1}{2} \left( \s_i - \s_j \s_k \right)^2- 1\,.
\eeq
followed by Hubbard-Stratonovich transformation. This leads to a non-interacting Hamiltonian ("$\s_i$"), and a two body interaction ("$\s_j \s_k$"), but one simply repeats the procedure on the latter, iterating until full decoupling is reached.

As a matter of fact, exact expressions for nonlinear models such as the perceptron are also possible: one Taylor expands the activation function, and applies the above procedure to the Taylor terms, one by one. }

\section*{Appendix}

Goetz Kersting\footnote{private communication.} points out that in case of two body interactions such as the SK-model, the Gaussian representation \eqref{exact_fe} is not in the least unique. The rationale is extremely charming, and best laid out by "reversed engineering": consider a centered Gaussian field $\{{\sf Z}_i \}_{i=1}^N$ with $\EE {\sf Z}_i {\sf Z}_j = g_{ij}, i \neq j$; the diagonal will be specified {\it a posteriori}. Undoing the $\log \cosh$, interchanging the order of integration, and by Gaussian integration we get
\beq \bea \label{gk}
\EE \exp \sum_i \log \cosh\left( h + \sqrt{\be} {\sf Z}_i \right) & = \EE E_o \exp\left( h \sum_i  \s_i + \sqrt{\be} \sum_i \s_i {\sf Z}_i \right) \\
&  = E_o \exp\left( h \sum_i \s_i+ \frac{\be}{2} \text{\sf var}\Big[\sum_i \s_i {\sf Z}_i \Big] \right) \,.
\eea \eeq
By definition of the $\sf Z$-field, and using that $\s_i^2 = 1$, it holds
\beq 
\text{\sf var}\left[\sum_i \s_i {\sf Z}_i \right] = {\sf E}\left[ \sum_{i,j} \s_i \s_j  {\sf Z}_i {\sf Z}_j \right] = \sum_{i \neq j} \s_i \s_j g_{ij} + \sum_i \EE {\sf Z}_i^2 \,.
\eeq
Plugging this in \eqref{gk} yields
\beq \bea
& \EE \exp\left( \sum_{i=1}^N \log \cosh\left( h + \sqrt{\be} {\sf Z}_i \right)  \right) = \\
& \hspace{3cm} = \exp \left( \frac{\be}{2} \sum_i \EE {\sf Z}_i^2 \right) \times E_o \exp\left( h \sum_i \s_i + \frac{\be}{2} \sum_{i \neq j} g_{ij} \s_i \s_j \right)\,.
\eea \eeq
The rightmost term above is the quenched partition function. Rearranging, this leads to the following representation of the quenched free energy for the SK-model:
\beq
N f_N(\be, h) = \log \EE \exp\left( \sum_{i=1}^N \log \cosh\left( h + \sqrt{\be} {\sf Z}_i \right)  \right) - \frac{\be}{2} \sum_{i=1}^N \EE {\sf Z}_i^2 \,.
\eeq
There is thus has a great deal of freedom  for the specification of the ${\sf Z}$-field: the only requirement is that the covariance matrix be positive semi-definite. Under this light, a most natural choice is a covariance given by ${\bf G} - \lambda_1 \mathbf 1 $, with $\la_1$ the smallest eigenvalue of the (Wigner) matrix $\bf G$. Consequences of this point of view will be worked out elsewhere.


\begin{thebibliography}{1} 

\bibitem{d} Domb, C., A. J. Barrett, and M. Lax. {\it Self-avoiding walks and real polymer chains.} Journal of Physics A: Mathematical, Nuclear and General 6.7 (1973): L82. 
\bibitem{ea} Edwards, Samuel Frederick, and Phil W. Anderson. {\it Theory of spin glasses.} Journal of Physics F: Metal Physics 5.5 (1975): 965.
\bibitem{fh} D. S. Fisher and D. A. Huse. {\it Absence of many states in realistic spin glasses}. J. Phys. A,
20(15):L1005–10, 1987.
\bibitem{gm} Gross, David J., and Marc M\'{e}zard. {\it The simplest spin glass.} Nuclear Physics B 240.4 (1984): 431-452.
\bibitem{he} Heisenberg, Werner. {\it Zur theorie des ferromagnetismus.} Original Scientific Papers Wissenschaftliche Originalarbeiten. Springer, Berlin, Heidelberg, 1985. 580-597.
\bibitem{hopfield} Hopfield, J. J. {\it Neural networks and physical systems with emergent collective computational abilities.} Proceedings of the National Academy of Sciences. 79 (8): 2554–2558. (1982)
\bibitem{h} Hubbard, John. {\it Calculation of partition functions.} Physical Review Letters 3.2 (1959): 77.
\bibitem{kac} M. Kac, {\it Trondheim Theoretical Physics Seminar}, Nordita Publ. No. 286, 1968 (unpublished); and T.-F. Lin, J. Math. Phys. 11, 1584 (1970).
\bibitem{mpv} Mézard, Marc, Giorgio Parisi, and Miguel A. Virasoro. {\it Spin glass theory and beyond.} World Scientific, Singapore (1987).
\bibitem{m} Montanari, Andrea. {\it Statistical mechanics and algorithms on sparse and random graphs.} Lectures on Probability Theory and Statistics. Saint-Flour (2013).
\bibitem{rosenblatt} Rosenblatt, Frank (1957). {\it The Perceptron—a perceiving and recognizing automaton}. Report 85-460-1. Cornell Aeronautical Laboratory.
\bibitem{sk} Sherrington, David, and Scott Kirkpatrick. {\it Solvable model of a spin-glass.} Physical review letters 35.26: 1792 (1975).


\end{thebibliography}
\end{document}